\documentclass[11pt]{amsart}
\usepackage{mathrsfs,amsthm,amssymb,amsmath}
\normalparindent=\parindent

\theoremstyle{plain}
    \newtheorem{thm}{Theorem}
    
    \newtheorem{lem}[thm]{Lemma}
    
    \newtheorem{cor}[thm]{Corollary}
    \newtheorem{fact}[thm]{Fact}
\theoremstyle{definition}
    \newtheorem{defn}[thm]{Definition}
    \newtheorem{nota}[thm]{Notation}
\theoremstyle{remark}
    \newtheorem{rem}[thm]{Remark}

\newcommand{\nin}{\notin}
\newcommand{\To}{\rightarrow}
\newcommand{\gdw}{\leftrightarrow}

 \DeclareMathOperator{\ma}{max}
\DeclareMathOperator{\med}{med} 

\newcommand{\cl}[1]{\langle #1 \rangle}
\newcommand{\clf}[1]{\langle \{ #1 \}\rangle}
\newcommand{\clfto}[1]{\langle \{ #1 \}\cup T_1 \rangle}

\newcommand{\Pow}{{\mathscr P}}
\renewcommand{\P}{{\mathscr P}}

\newcommand{\C}{{\mathscr C}}
\newcommand{\R}{{\mathscr R}}
\newcommand{\F}{{\mathscr F}}
\newcommand{\J}{{\mathscr J}}
\newcommand{\A}{{\mathscr A}}
\newcommand{\B}{{\mathscr B}}
\newcommand{\D}{{\mathscr D}}
\newcommand{\M}{{\mathscr M}}
\newcommand{\N}{{\mathscr N}}
\newcommand{\U}{{\mathscr U}}

\newcommand{\OO}{{\mathscr O}}
\renewcommand{\O}{{\mathscr O}}

\newcommand{\On}{{\mathscr O}^{(n)}}

\newcommand{\Oo}{{\mathscr O}^{(1)}}
\newcommand{\Ot}{{\mathscr O}^{(2)}}

\newcommand{\cto}{\langle T_1 \rangle}
\newcommand{\pto}{Pol(T_1)}
\newcommand{\pton}{\pto^{(n)}}

\newcommand{\niceint}{[\U,\O]}
\newcommand{\lw}{\leq_W}
\newcommand{\aw}{\sim_{W}}
\newcommand{\fto}[1]{\langle #1\rangle_{T_1}}
\newcommand{\un}{^{(n)}}
\newcommand{\uk}{^{(k)}}
\newcommand{\ut}{^{(2)}}

\newcommand{\ow}{\text{otherwise}}

\author[M.\,Pinsker]{Michael Pinsker}
\address{Algebra\\TU Wien\\Wiedner Hauptstrasse 8-10/118\\A-1040 Wien, Austria}
\email{marula@gmx.at}
\title[Almost unary functions]{Clones containing all almost unary functions}
\subjclass{Primary 08A40; secondary 08A05}

\keywords{clone lattice, maximal clones, polymorphisms, median
functions, almost unary functions, regular cardinals}

\thanks{The author is supported by DOC [Doctoral Scholarship Programme of the Austrian Academy of Sciences].
He is grateful to M. Goldstern for drawing his attention to the
subject and for his remarks on the paper, to L. Heindorf for his
comments and to the II. Mathematisches Institut at Freie
Universit\"{a}t Berlin for their hospitality during his visit}

\begin{document}

\begin{abstract}
    Let $X$ be an infinite set of regular cardinality. We determine all clones on $X$ which contain all almost
    unary functions. It turns out that independently of the size of $X$, these clones form a countably infinite descending
    chain. Moreover, all such clones are finitely generated over the unary functions. In particular,
    we obtain an explicit description of the only maximal clone in this part of the clone lattice. This is especially
    interesting if $X$ is countably infinite, in which case it is known that such a description cannot be obtained
    for the second maximal clone over the unary functions.
\end{abstract}
\maketitle \thispagestyle{empty}
\setlength{\baselineskip}{1.1\baselineskip}

\begin{section}{Introduction}

\subsection{Clones and almost unary functions}
    Let $X$ be a set and denote by $\On$ the set of all $n$-ary functions
    on $X$. Then $\OO=\bigcup_{n=1}^{\infty}\On$ is the set of all functions
    on $X$. A \emph{clone} $\C$ over $X$ is a subset of $\OO$ which
    contains the projections and which is closed under compositions.
    The set of all clones over $X$ forms a complete lattice $Clone(X)$ with
    respect to inclusion. This lattice is a subset of the power set of $\OO$.
    The clone lattice is countably infinite if
    $X$ has only two elements, but is of size $2^{\aleph_0}$ for
    cardinality of $X$ finite and greater than two. For infinite $X$ we have
    $|Clone(X)|=2^{2^{|X|}}$.

    Let $X$ be of infinite regular cardinality from now on unless otherwise stated. We
    call a subset $S\subseteq X$ \emph{large} iff $|S|=|X|$, and
    \emph{small} otherwise. If $X$ is itself a regular cardinal,
    then the small subsets are exactly the bounded subsets of $X$.
    A function $f(x_1,...,x_n)\in\On$ is \emph{almost unary} iff
    there exists a function $F:X\To \Pow(X)$ and $1\leq k\leq n$ such that $F(x)$ is small for all
    $x\in X$ and such that for all
    $(x_1,...,x_n)\in X^n$ we have $f(x_1,...,x_n)\in F(x_k)$. If
    we assume $X$ to be a regular cardinal itself, this is
    equivalent to the existence of a function $F\in\Oo$ and a
    $1\leq k\leq n$ such that $f(x_1,...,x_n)< F(x_k)$ for all $(x_1,...,x_n)\in
    X^n$. Because this is much more convenient and does not
    influence the properties of the clone lattice, we shall assume
    $X$ to be a regular cardinal throughout this paper. Let $\U$ be the set of all almost unary functions. It is
    readily verified that $\U$ is a clone (see e.g. \cite{DR85}). We will determine all clones which contain $\U$; in particular, such clones contain $\Oo$.

    \subsection{Maximal clones above $\Oo$}
    A clone is called \emph{maximal}
    iff it is a dual atom in $Clone(X)$. For finite $X$, the number of maximal clones is
    finite and all of them are known (a result due to I. Rosenberg
    \cite{Ros70}, see also \cite{Pin02}). Moreover, the lattice is dually atomic in that case,
    i.e. every clone is contained in a maximal one. If $X$ is infinite, the
    situation is rather hopeless as another theorem by I. Rosenberg
    \cite{Ros76} states that there exist $2^{2^{|X|}}$ maximal
    clones, see also \cite{GS022}.
    In addition, a recent result due to M. Goldstern and S. Shelah
    \cite{GS03} shows that if the continuum hypothesis holds, then the clone
    lattice of a countable base set is not even dually atomic.

    However, by Zorn's lemma, the sublattice of $Clone(X)$ of functions containing $\Oo$ is dually
    atomic since $\OO$ is finitely generated over $\Oo$.
    G. Gavrilov proved in \cite{Gav65} that for countably
    infinite $X$ there are only two maximal clones containing all unary
    functions. M. Goldstern and S. Shelah extended this result to
    clones on weakly compact cardinals in the article \cite{GS022}, but proved
    also that on other regular cardinals $X$ satisfying a certain partition
    relation there are even $2^{2^X}$ such clones.

    There exists exactly one maximal clone above $\U$. So far,
    this clone has been defined using the following concept:
    Let $\rho\subseteq X^J$ be a relation on $X$ indexed by $J$ and let
    $f\in\On$. We say that $f$ \emph{preserves} $\rho$ iff for all $r^1=(r^1_i:i\in J),\cdots,r^n=(r^n_i:i\in J)$ in
    $\rho$ we have $(f(r^1_i,...,r^n_{i}):i\in J) \in \rho$.
    For a set of relations $\R$ on $X$ we define the set of
    \emph{polymorphisms} $Pol(\R)$ of $\R$ to be the set of all
    functions in $\OO$ preserving all $\rho\in\R$. In particular, if
    $\rho\subseteq X^{X^k}$ is a set of $k$-ary functions, then a
    function $f\in\On$ preserves $\rho$ iff for all functions
    $g_1,...,g_n$ in $\rho$ the composite $f(g_1,...,g_n)$ is a
    function in $\rho$.

    Write
    $$
        T_1=\U\ut=\{f \in \Ot: f \,\,\text{almost unary}\}.
    $$

    The following was observed by G. Gavrilov \cite{Gav65} for
    countable base sets and extended to all regular $X$ by R. Davies and I.
    Rosenberg \cite{DR85}. Uniqueness on uncountable regular cardinals
    is due to M. Goldstern and S. Shelah \cite{GS022}.
    \begin{fact}\label{FAC:poltone}
        Let $X$ have infinite regular cardinality. Then $Pol(T_1)$ is a maximal clone containing all unary
        functions. Furthermore, $\pto$ is the only maximal clone
        containing all almost unary functions.
    \end{fact}

    For $S$ a subset of $X$ we set
    $$
        \Delta_S=\{(x,y)\in S^2 : y < x\},\quad \nabla_S=\{(x,y)\in S^2 : x
        <y\}.
    $$
    We will also write $\Delta$ and $\nabla$ instead of $\Delta_X$ and $\nabla_X$. Now define
    $$
        T_2=\{f\in\Ot : \forall S\subseteq X \,(S \,\,\text{large}
        \To\,
        \text{neither} \, f\upharpoonright_{\Delta_S}\,\text{nor} \,
        f\upharpoonright_{\nabla_S}\, \text{are 1-1})\}.
    $$
    The next result is due to G. Gavrilov \cite{Gav65} for $X$ a countable set and due to M. Goldstern and S. Shelah
    \cite{GS022} for $X$ weakly compact.
    \begin{fact}\label{FAC:theTwoMaximal}
         Let $X$ be countably infinite or weakly compact.
         Then $Pol(T_2)$ is a maximal clone which contains $\Oo$. Moreover,
         $Pol(T_1), Pol(T_2)$ are the only maximal clones above $\Oo$.
    \end{fact}

    The definition of $Pol(T_2)$ not only looks more
    complicated than the one of $Pol(T_1)$. First of all, a result of R. Davies and I.
    Rosenberg in \cite{DR85} shows that assuming the continuum hypothesis, $T_2$ is not closed under
    composition on $X=\aleph_1$ and so it is unclear what $Pol(T_2)$ is.
    Secondly, on countable $X$, if we equip $\OO$ with a certain natural topology which we shall specify later,
    then $T_2$ is a complete $\Pi_1^1$ set in that space and so is
    $Pol(T_2)$; in particular, neither $T_2$ nor $Pol(T_2)$ are countably generated over
    $\Oo$ (see \cite{Gol02}). The set $T_1$ on the other hand is generated by a single binary function over $\Oo$: Let $p$ be
    any injection from $X^2$ to $X$. For technical reasons we assume that $0$ is not in the range of
    $p$. The next fact is folklore.
    \begin{fact}
        $\cl{\{p\}\cup\Oo}=\O$.
    \end{fact}
    For a subset
    $S$ of $X^2$ we write
    $$
        p_S(x_1,x_2)=\begin{cases}p(x_1,x_2)&,(x_1,x_2)\in
        S\\0&,\ow\\ \end{cases}
    $$
    M. Goldstern observed the following \cite{Gol02}. Since the
    result has not yet been published, but is important for our investigations, we include a proof here.
    \begin{fact}\label{FAC:pDelta}
        $\cl{\{p_\Delta\}\cup\Oo}=\cto$.
    \end{fact}
    \begin{proof}
        Set $\C=\cl{\{p_\Delta\}\cup\Oo}$. Since $p_\Delta(x_1,x_2)$ is
        obviously bounded by the unary function
        $\gamma(x_1)=\sup\{p_\Delta(x_1,x_2):x_2\in
        X\}+1=\sup\{p(x_1,x_2):x_2<x_1\}+1$, where by $\alpha+1$ we
        mean the successor of an ordinal $\alpha$,
        we have $p_\Delta\in T_1$ and hence
        $\C\subseteq\cto$.

        To see the other inclusion, note first that the function
        $$
            q(x_1,x_2) =\begin{cases}p_\Delta(x_1,x_2)&, (x_1,x_2)\in\Delta\\
                 x_1&,\ow\end{cases}
        $$
        is in $\C$. Indeed, choose $\epsilon\in\Oo$ strictly
        increasing such that $p_\Delta(x_1,x_2)<\epsilon(x_1)$ for
        all $x_1,x_2\in X$ and consider
        $t(x_1,x_2)=p_\Delta(\epsilon(x_1),p_\Delta(x_1,x_2))$. On
        $\Delta$, $t$ is still one-one, and outside $\Delta$, the term
        is a one-one function of the first component $x_1$. Moreover, the ranges $t[\Delta]$ and
        $t[X^2\setminus\Delta]$ are disjoint. Hence,
        we can write $q=u\circ t$ for some unary $u$.
        By the same argument we see that for arbitrary unary functions $a,b\in\Oo$ the
        function
        $$
                q_{a,b}(x_1,x_2) =\begin{cases}a(p_\Delta(x_1,x_2))&, (x_1,x_2)\in\Delta\\
                 b(x_1)&,\ow\end{cases}
        $$
        is an element of $\C$.

        Now let $f\in T_1$ be given and say $f(x_1,x_2) < \delta(x_1)$ for all $x_1,x_2\in X$,
        where $\delta\in\Oo$ is strictly increasing. Choose $a\in\Oo$ such that $a(p_\Delta(x_1,x_2))=f(x_1,x_2)+1$ for all
        $(x_1,x_2)\in \Delta$. Then set
        $$
        f_1(x_1,x_2) =q_{a,\delta+1}(x_1,x_2)=\begin{cases}
              f(x_1,x_2)+1&,(x_1,x_2)\in\Delta\\
                \delta(x_1)+1           &,\ow
                \end{cases}
        $$
        We construct a second function
        $$
        f_2(x_1,x_2) =\begin{cases} 0&, (x_1,x_2)\in\Delta\\
                                f(x_1,x_2)+1&, \ow\end{cases}
        $$
        It is readily verified that
        $f_2(x_1,x_2) = u (p_\Delta(x_2+1, x_1))$ for some unary $u$.
        Now  $f_2(x_1,x_2) <  f_1(x_1,x_2)$ and $f_1,f_2\in\C$.
        Clearly $$f(x_1,x_2)  = u(p_\Delta( f_1(x_1,x_2) , f_2(x_1,x_2) ))$$  for some
        unary $u$. This shows $f\in\C$ and so $\cto\subseteq\C$ as
        $f\in T_1$ was arbitrary.
    \end{proof}

    We shall see that $\pto$ is also finitely generated
    over $\Oo$. Moreover, for countable $X$ it is a Borel set in the topology yet to be
    defined.
    Our explicit description $\pto$ holds for all infinite $X$ of regular cardinality, but
    is interesting only if there are not too many other maximal
    clones containing $\Oo$. By Fact \ref{FAC:theTwoMaximal}, this
    is at least the case for $X$ countably infinite or weakly
    compact.

\subsection{Notation}

    For a set of functions $\F$ we shall denote the smallest
    clone containing $\F$ by $\langle \F \rangle$. By $\F^{(n)}$ we refer to the set
    of $n$-ary functions in $\F$. \\
    We call the projections which every clone contains $\pi^n_i$
    where $n\geq 1$ and $1\leq i\leq n$.
    If $f\in\On$ is an $n$-ary function, it sends $n$-tuples of elements of $X$ to $X$ and we
    write $(x_1,...,x_n)$ for these tuples unless otherwise stated as in $f(x,y,z)$; this is the only place where
    we do not stick to set-theoretical notation (according to which we would have to write $(x_0,...,x_{n-1})$).
    The set $\{1,...,n\}$ of indices of $n$-tuples will play an important role and we write $N$
    for it. We denote the set-theoretical complement of
    a subset $A\subseteq N$ in $N$ by $-A$. We identify the set $X^n$ of $n$-tuples
    with the set of functions from $N$ to $X$, so that if $A\subseteq N$ and $a: A\To X$ and $b: -A\To X$ are
    partial functions, then $a\cup b$ is an $n$-tuple. Sometimes, if the arity of $f\in\OO$ has not yet been given a name,
    we refer to that arity by $n_f$.\\
    If $a\in X^n$ is an $n$-tuple and $1\leq k\leq n$ we write
    $(a)^n_k$ or only $a_k$ for the $k$-th component of $a$. For $c\in
    X$ and $J$ an index set we write $c^J$ for the $J$-tuple with
    constant value $c$. The
    order relation $\leq$ on $X$ induces the
    pointwise partial order on the set of $J$-tuples of elements of $X$ for any index
    set $J$: For $x,y\in X^J$ we write $x\leq y$ iff $x_j\leq y_j$ for all $j\in J$. Consequently we also
    denote the induced pointwise
    partial order of $\On$ by $\leq$, so that for
    $f,g\in\On$ we have $f\leq g$ iff $f(x)\leq g(x)$ for all $x \in
    X^n$. Whenever we state that a function $f\in\On$ is monotone, we
    mean it is monotone with respect to $\leq$: $f(x)\leq f(y)$
    whenever $x\leq y$. We denote the power set of $X$ by $\Pow(X)$. The
    element $0\in X$ is the smallest element of $X$.
\end{section}

\begin{section}{Properties of clones above $\U$ and the clone $\pto$}

\subsection{What $\cto$ is}
    We start by proving that the almost unary clone $\U$ is a so-called \emph{binary
    clone}, that is, it is generated by its binary part. Thus,
    when investigating $[\U,Pol(T_1)]$, we are in fact dealing
    with an interval of the form $[\cl{\C^{(2)}},Pol(\C^{(2)})]$ for $\C$ a clone.
    \begin{lem}
        The binary almost unary functions generate all almost unary
        functions. That is, $\cto=\U$.
    \end{lem}
    \begin{proof}
        Trivially, $\cto\subseteq \U$. Now we prove by induction
        that $\U^{(n)}\subseteq \cto$ for all $n\geq 1$. This is
        obvious for $n=1,2$. Assume we have $\U^{(k)}\subseteq
        \cto$ for all $k<n$ and take any function $f\in\U^{(n)}$.
        Say without loss of generality that $f(x_1,...,x_n)\leq\gamma
        (x_1)$ for some $\gamma\in\Oo$. We will use the function
        $p_\Delta\in T_1$ to code two variables into one and then
        use the induction hypothesis. Define
        $$
            g_1(x_1,...,x_{n-2},z)=
            \begin{cases}f(x_1,...,x_{n-2},(p_\Delta^{-1}(z))^2_1,(p_\Delta^{-1}(z))^2_2)\quad
            &,z\in p_\Delta[X^2]\setminus\{0\} \\
            x_1\quad &,\ow \end{cases}
        $$
        The function is an element of $\U^{(n-1)}$ as it is bounded by
        $\max(x_1,\gamma(x_1))$. Intuitively, $g_1$ does the following: If $z\neq 0$ and in the range
        of $p_\Delta$, then
        $g_1$ imagines a pair $(x_{n-1},x_n)$ to be coded into $z$
        via $p_\Delta$. It reconstructs the pair $(x_{n-1},x_n)$
        and calculates $f(x_1,...,x_n)$. If $z=0$ or not in the range of $p_\Delta$, then $g$ knows
        there is no information in $z$; it simply forgets about the tuple
        $(x_2,...,x_n)$ and returns $x_1$, relying on the following similar function
        to do the job: Set $\Delta'=\Delta\cup\{(x,x):x\in X\}$
        and define
        $$
            g_2(x_1,...,x_{n-2},z)=
            \begin{cases}f(x_1,...,x_{n-2},(p_{\Delta'}^{-1}(z))^2_2,(p_{\Delta'}^{-1}(z))^2_1)\quad
            &,z\in p_{\Delta'}[X^2]\setminus\{0\}\\
            x_1\quad &,\ow\end{cases}
        $$
        The function $g_2$ does exactly the same as $g_1$ but assumes the
        pair $(x_{n-1},x_n)$ to be coded into $z$ in wrong
        order, namely as $(x_n,x_{n-1})$, plus it cares for the diagonal.
        Now consider
        $$
            h(x_1,...,x_n)=g_2(g_1(x_1,...,x_{n-2},p_\Delta(x_{n-1},x_n)),x_2,...,x_{n-2},p_{\Delta'}(x_n,x_{n-1})).
        $$
        All functions which occur in $h$ are almost unary with at most $n-1$ variables.
        We claim that $h=f$. Indeed, if $x_{n-1}<x_n$, then
        $p_\Delta(x_{n-1},x_n)\neq 0$ and $g_1$ yields $f$. But
        $p_{\Delta'}(x_n,x_{n-1})=0$ and so $g_2$ returns $g_1=f$. If
        on the other hand $x_n\leq x_{n-1}$, then $p_\Delta(x_{n-1},x_n)=
        0$ and $g_1=x_1$, whereas $p_{\Delta'}(x_n,x_{n-1})\neq 0$,
        which implies $g_2=f(g_1,x_2,...,x_n)=f(x_1,...,x_n)$.
    \end{proof}
    The following lemma will be
    crucial for our investigation of clones containing $T_1$.
    \begin{cor}
        Let $\C$ be a clone containing $T_1$. Then $\C$ is
        downward closed, that is, if $f\in \C$, then also $g\in \C$
        for all $g\leq f$.
    \end{cor}
    \begin{proof}
        If $f\in \C^{(n)}$ and $g\in\On$ with $g\leq f$ are given, define
        $h_g(x_1,...,x_{n+1})=\min(g(x_1,...,x_n),x_{n+1})$. Then
        $h_g\leq x_{n+1}$ and consequently, $h_g\in\cto\subseteq\C$.
        Now $g=h_g(x_1,...,x_n,f(x_1,...,x_n))\in\C$.
    \end{proof}

\subsection{Wildness of functions}
    We have seen in the last section that the interval $\niceint$ is
    about growth of functions as all clones in that interval
    are downward closed. But mind we are not talking about how
    rapidly functions are growing in the sense of polynomial growth,
    exponential growth and so forth since we are considering
    clones modulo $\Oo$ (and so we can make functions as steep as
    we like); the growth of a function will be determined by which
    of its variables are responsible for the function to obtain many values.
    The following definition is due to M. Goldstern and S. Shelah
    \cite{GS022}. Recall that $N = \{1,...,n\}$.
    \begin{defn}
        Let $f\in\On$. We call a set $\varnothing\neq A\subseteq N$ \emph{$f$-strong} iff for
        all $a\in X^A$ the set $\{f(a\cup x):x\in X^{-A}\}$ is
        small. $A$ is \emph{$f$-weak} iff it is not $f$-strong.
        In order to use the defined notions more freely, we define the
        empty set to be $f$-strong iff $f$ has small range.
    \end{defn}
    Thus, a set of indices of variables of $f$ is strong iff $f$ is
    bounded whenever those variables are. For example, a function is almost unary
    iff
    it has a one-element strong set. Here, we shall rather
    think in terms of the complements of weak sets.
    \begin{defn}
        Let $f\in\On$ and let $A\subsetneqq N$ and $a\in X^{-A}$. We
        say $A$ is \emph{$(f,a)$-wild} iff the set $\{f(a\cup
        x):x\in X^A\}$ is large. The set $A$ is called
        \emph{$f$-wild} iff there exists $a\in X^{-A}$ such that $A$ is
        $(f,a)$-wild. We say that $A$ is \emph{$f$-insane} iff $A$
        is $(f,a)$-wild for all $a\in X^{-A}$. The set $N$ itself
        we call $f$-wild and $f$-insane iff $f$ is unbounded.
    \end{defn}
    Observe that if $A\subseteq B\subseteq N$ and $A$ is $f$-wild, then $B$ is $f$-wild as well.
    Obviously, $A\subseteq N$ is $f$-wild iff $-A$ is
    $f$-weak. It is useful to state the following trivial criterion for a
    function to be almost unary.
    \begin{lem}\label{LEM:almostUnaryCriterion}
        Let $n\geq 2$ and $f\in\On$. $f$ is almost
        unary iff there exists a subset of $N$ with $n-1$ elements which is not $f$-wild.
    \end{lem}
    \begin{proof}
        If $f$ is almost unary, then there is a one-element $f$-strong subset of $N$ and
        the complement of that set is not $f$-wild. If on the other hand there exists
        $k\in N$ such that $N\setminus\{k\}$ is
        not $f$-wild, then $\{k\}$ is $f$-strong and so $f$ is almost unary.
    \end{proof}

    We will require the following fact from
    \cite{GS022}.
    \begin{fact} \label{FAC:goldsternShelah}
        If $f\in\pton$ and $A_1,A_2\subseteq N$ are $f$-wild, then
        $A_1\cap A_2\neq \varnothing$.
    \end{fact}
    We observe that the converse of this statement holds as well.
    \begin{lem}\label{charPto}
        Let $f\in\On$ be any $n$-ary function. If all pairs of
        $f$-wild subsets of $N$ have a nonempty intersection, then
        $f\in\pto$.
    \end{lem}
    \begin{proof}
        Let $g_1,...,g_n\in T_1$ be given and set $A_1=\{k\in
        N:\exists\gamma\in\Oo\,(g_k(x_1,x_2)\leq\gamma(x_1))\}$ and
        $A_2=-A_1$. Since $A_1\cap A_2=\varnothing$ either $A_1$
        or $A_2$ cannot be $f$-wild. Thus
        $f(g_1,...,g_n)$ is bounded by a unary function of $x_2$
        in the first case and by a unary function of $x_1$ in the
        second case.
    \end{proof}
    The equivalence yields a first description of
    $\pto$ with an interesting consequence.
    \begin{thm}\label{COR:FirstptoDescription}
        A function $f\in\On$ is an element of $\pto$ iff all pairs
        of $f$-wild subsets of $N$ have a nonempty intersection.
    \end{thm}

    We show now that for countable $X$, this description implies that $\pto$ is a Borel set with respect to
    the natural topology on
    $\OO$. We do not explain the basic
    notions of descriptive set theory. The reader not familiar
    with these notions is advised to either skip this part and proceed directly to
    the next section or to consult \cite{Kec95}.

    Equip $X=\omega$ with the
    discrete topology. Then the product space $\N=\omega^\omega=\Oo$
    is the so-called Baire space. It is obvious that
    $\On=\omega^{\omega^n}$ is homeomorphic to $\N$. Thus,
    $\OO=\bigcup_{n=1}^{\infty}\On$ is the sum of $\omega$ copies of
    $\N$.

    \begin{thm}
        Let $X$ be countably infinite. Then $\pto$ is a Borel set
        in $\OO$.
    \end{thm}
    \begin{proof}
        By the preceding theorem,
        $$
            \pto\un=\{f\in\On: \forall A,B\subseteq N
            (A,B \,f\text{-wild} \To A\cap B\neq \varnothing )\}
        $$
        There are no (only finite) quantifiers in this definition
        except for those which might occur in the predicate of wildness.
        Now
        $$
            A\subseteq N\,f\text{-wild} \gdw \exists a\in X^{-A}\forall
            k\in X\exists b\in X^A (f(a\cup b)>k)
        $$
        All quantifiers range over countable sets so that $\pto$ is
        $\Sigma_3^0$.
    \end{proof}
    The preceding theorem shows that $Pol(T_2)$ is much more
    complicated than $Pol(T_1)$, as M. Goldstern observed the
    following \cite{Gol02}.
    \begin{fact}
        Let $X$ be countably infinite. Then $Pol(T_2)$ is a $\Pi_1^1$-complete set in $\O$.
    \end{fact}

\subsection{What wildness means}
    We wish to compare the wildness of functions. Write $S_N$ for the set of all permutations on $N$.
    \begin{defn}
        For $f,g\in\On$ we say that $f$ is \emph{as wild as $g$}
        and write $f \aw g$ iff there exists a permutation $\pi\in
        S_N$ such that $A$ is $f$-wild if and only
        if $\pi[A]$ is $g$-wild for all $A\subseteq N$. Moreover, $g$ is \emph{at least as wild
        as $f$} ($f\lw g$) iff there is a permutation $\pi\in S_N$ such that
        for all $f$-wild subsets $A\subseteq N$ the image $\pi[A]$
        of $A$ under $\pi$ is $g$-wild.
    \end{defn}

    \begin{lem}
        $\aw$ is an equivalence relation and $\lw$ a quasiorder extending $\leq$ on
        the set of $n$-ary functions $\On$.
    \end{lem}
    \begin{proof}
        We leave the verification of this to the reader.
    \end{proof}
    \begin{lem}
        Let $f,g\in\On$. Then $f\aw g$ iff $f\lw g$ and $g\lw f$.
    \end{lem}
    \begin{proof}
        It is clear that $f\lw g$ (and $g\lw f$) if $f\aw g$. Now
        assume $f\lw g$ and $g\lw f$. Then there are
        $\pi_1,\pi_2\in S_N$ which take $f$-wild and $g$-wild
        subsets of $N$ to $g$-wild and $f$-wild sets, respectively.

        Set $\pi=\pi_2\circ\pi_1$. Then $A$ is $f$-wild iff $\pi[A]$
        is $f$-wild for any subset $A$ of $N$: If $A$ is $f$-wild,
        then $\pi_1[A]$ is $g$-wild, then $\pi_2[\pi_1[A]]=\pi[A]$
        is $f$-wild. If on the other hand $\pi[A]$ is $f$-wild,
        then take $k\geq 1$ such that $\pi^k=id_N$ and observe
        that $\pi^{k-1}\circ\pi[A]=\pi^k[A]=A$ is $f$-wild.

        Now we see that $A$ is $f$-wild iff $\pi_1[A]$ is $g$-wild for all $A\subseteq N$:
        If $\pi_1[A]$ is $g$-wild, then so is
        $\pi_2\circ\pi_1[A]=\pi[A]$ and so is $A$ by the preceding
        observation. Hence, the permutation $\pi_1$ shows that
        $f\aw g$.
    \end{proof}

    \begin{cor}
        Let $n\geq 1$. Then $\mathord{\leq_W}/\mathord{\aw}$ is a partial order on
        the $\aw$-equivalence classes of $\On$.
    \end{cor}

    \begin{nota}
        Let $f \in\On$. By $\fto{f}$ we mean $\clfto{f}$ from now
        on. $\fto{f}$ is the smallest clone containing $f$ as well
        as all almost unary functions.
    \end{nota}
    We are aiming for the following theorem which tells us why we invented wildness.
    \begin{thm}\label{THM:mainWild}
        Let $f,g\in \On$. If $f\lw g$, then $f\in \fto{g}$. In
        words, if $g$ is at least as wild as $f$, then it
        generates $f$ modulo $T_1$.
    \end{thm}
    \begin{cor}
        Let $f,g\in\On$. If $f\aw g$, then $\fto{f}=\fto{g}$.
    \end{cor}
    We split the proof of Theorem \ref{THM:mainWild} into a
    sequence of lemmas. In the next lemma we see that it does not matter which $a\in
    X^{-A}$ makes a set $A\subseteq N$ wild.
    \begin{lem}\label{zerowild}
        Let $g\in\On$. Then there exists $g'\in\fto{g}\un$ such
        that for all $A\subseteq N$ the following holds:  If $A$ is $g$-wild, then $A$
        is $(g',0^{-A})$-wild.
    \end{lem}
    \begin{proof}
        Fix for all $g$-wild $A\subseteq N$ a tuple $a_A\in
        X^{-A}$ such that $\{g(x\cup a_A):x\in X^A\}$ is large.
        For an $n$-tuple $(x_1,...,x_n)$ write $P=P(x_1,...,x_n)=\{l\in N:x_l\neq
        0\}$ for the set of indices of positive components in the tuple.
        Define for $1\leq i\leq n$ functions
        $$
            \gamma_i(x_1,...,x_n)=\begin{cases}x_i&,x_i\neq 0 \vee
            P(x_1,...,x_n)\,\,\text{not}\,\,
            g\text{-wild}\\(a_{P})_i&,\text{otherwise}\end{cases}
        $$
        In words, if the set $P$ of indices of positive components in $(x_1,...,x_n)$ is a wild set,
        then the $\gamma_i$ leave those positive components alone and send the zero components
        to the respective values making $P$ wild. Otherwise, they act just like projections.
        It is obvious that $\gamma_i$ is almost unary, $1\leq i \leq
        n$. Set $g'=g(\gamma_1,...,\gamma_n)\in\fto{g}$. To
        prove that $g'$ has the desired property, let $A\subseteq N$
        be $g$-wild. Choose any minimal $g$-wild $A'\subseteq A$. Then
        by the definition of wildness the set $\{g(x\cup a_{A'}):x\in
        X^{A'}\}$ is large. Take a large
        $B\subseteq X^{A'}$ such that the sequence $(g(x\cup
        a_{A'}): x\in B)$ is one-one. Select further a large $C\subseteq
        B$ such that each component in the sequence of tuples
        $(x:x\in C)$ is either constant or injective and such that $0$ does not occur
        in any of the injective components (it is a simple
        combinatorial fact that this is possible). If one of
        the components were constant, then $A'$
        would not be minimal $g$-wild; hence, all components are
        injective. Now we have
        $$
        \begin{aligned}
            |X|&=|\{g(x\cup a_{A'}): x\in C\}|\\&= |\{g'(x\cup
            0^{-A'}):x\in C\}|
            &\leq|\{g'(x\cup 0^{-A}):x\in X^A\}|
        \end{aligned}
        $$
        and so $A$ is $(g',0^{-A})$-wild.
    \end{proof}
    We prove that we can assume functions to be monotone.
    \begin{lem}\label{monotone}
        Let $g\in\On$. Then there exists $g''\in\fto{g}\un$ such
        that $g\leq g''$ and $g''$ is monotone with respect to the
        pointwise order $\leq$.
    \end{lem}
    \begin{proof}
    We will define a mapping $\gamma$ from $X^n$ to $X^n$ such that $\gamma_i=\pi^n_i\circ\gamma$ is almost unary
        for $1\leq i\leq n$ and such that
        $g''=g\circ\gamma$ has the desired property. We fix for every
        $g$-wild $A\subseteq N$ a sequence $(\alpha^A_\xi)_{\xi\in X}$
        of elements of $X^n$ so that all components of $\alpha^A_\xi$ which lie not in $A$ are constant and so that
        $(g(\alpha^A_\xi))_{\xi\in X}$ is
        monotone and unbounded.

        Let $x\in X^n$. The \textit{order type} of $x$ is the unique
        $n$-tuple $(j_1,...,j_n)$ of indices in $N$ such that
        $\{j_1,...,j_n\}=\{1,...,n\}$ and such that
        $x_{j_1}\leq...\leq x_{j_n}$ and such that $j_k<j_{k+1}$ whenever
        $x_{j_k}=x_{j_{k+1}}$. Let $1\leq k\leq n$ be the largest element with the property that the set
        $\{j_k,...,j_n\}$ is $g$-wild. We call the set $\{j_k,...,j_n\}$ the \textit{pushing
        set} $Push(x)$
        and $\{j_1,...,j_{k-1}\}$ the \textit{holding set} of $x$ with
        respect to $g$.

        We define by transfinite recursion
        $$
            \gamma:\quad\begin{matrix} X^n &\To& X^n\\
            x&\mapsto& \alpha^{Push(x)}_{\lambda(x)} \end{matrix}
        $$
        where
        $$
            \lambda(x)=\min\{\xi:g(\alpha^{Push(x)}_\xi)\geq\sup(\{g''(y):y<x\}\cup\{g(x)\})\}.
        $$
        This looks worse than it is: We simply map $x$ to the
        first element of the sequence $(\alpha^{Push(x)}_\xi)_{\xi\in X}$ such that all values of $g''$
        already defined as well as $g(x)$ are topped. By
        definition, $g''=g\circ\gamma$ is monotone and $g\leq
        g''$. It only remains to prove that all $\gamma_i$, $1\leq
        i\leq n$, are almost unary to see that $g''\in\fto{g}$.

        Suppose not, and say that $\gamma_k$ is not almost unary for some $1\leq k\leq n$. Then there exists a
        value $c\in X$ and a sequence of $n$-tuples
        $(\beta_\xi)_{\xi\in X}$ with constant value $c$ in the $k$-th component
        such that $(\gamma_k(\beta_\xi))_{\xi\in X}$ is unbounded.
        Since there exist only finitely many order types of
        $n$-tuples, we can assume that all $\beta_\xi$ have the same
        order type $(j_1,...,j_n)$; say without loss of generality $(j_1,...,j_n)=(1,...,n)$.
        Then all $\beta_\xi$ have the same pushing set $Push(\beta)$ of indices.
        If $k$ was an element of the holding set of the
        tuples $\beta_\xi$, then $(\gamma_k(\beta_\xi):\xi\in X)$ would be constant so that $k$ must be in $Push(\beta)$.
        Clearly, $(\lambda(\beta_\xi))_{\xi\in X}$
        has to be unbounded as otherwise $(\gamma_k(\beta_\xi))_{\xi\in X}$ would be bounded. Since by definition
        the value of $\lambda$
        increases only when it is necessary to keep $g\leq g''$,
        the set $\{g(y):\exists \xi\in X(y\leq \beta_\xi)\}$ is unbounded. But because of the order type of the $\beta_\xi$,
        whenever $i\leq k$, then we have $(\beta_\xi)^n_{i}\leq c$ for all
        $\xi\in X$ so that the components of the $\beta_\xi$ with index in the set
        $\{1,...,k\}$ are bounded. Thus, $\{k+1,...,n\}$
        is $g$-wild, contradicting the fact that $k$ is in the
        pushing set $Push(\beta)$.

    \end{proof}
    In a next step we shall see that modulo $T_1$, wildness is insanity.
    \begin{lem}\label{LEM:MonotoneAndInsane}
        Let $g\in\On$. Then there exists $g''\in\fto{g}\un$ such
        that $g''$ is monotone and for all $A\subseteq N$ the following holds:  If $A$ is $g$-wild, then $A$
        is $g''$-insane.
    \end{lem}
    \begin{proof}
        Let $g'\in\fto{g}\un$ be provided by Lemma \ref{zerowild} and make a monotone $g''$ out
        of it with the help of the preceding lemma. We claim that
        $g''$ already has both desired properties. To prove this,
        consider an arbitrary $g$-wild $A\subseteq N$. By construction
        of $g'$, $A$ is $(g',0^{-A})$-wild and so it is also
        $(g'',0^{-A})$-wild as $g'\leq g''$. But $0^{-A}\leq a$
        for all $a\in X^{-A}$; hence the fact that $g''$ is monotone implies
        that $A$ is $(g,a)$-wild for all $a\in X^{-A}$ which means
        exactly that $A$ is $g''$-insane.
    \end{proof}
    \begin{lem}\label{leq}
        Let $f, g\in\On$. If $f\lw g$, then there exists
        $h\in\fto{g}\un$ such that $f\leq h$.
    \end{lem}
    \begin{proof}
        Without loss of generality, we assume that the permutation $\pi\in S_N$ taking
        $f$-wild subsets of $N$ to $g$-wild sets is the identity on $N$.
        We take $g''\in\fto{g}$ according to the preceding lemma.
        We wish to define $\gamma\in\Oo$ with $f\leq\gamma\circ
        g''$. For $x\in X$ write $U_x=g''^{-1}[\{x\}]$ for the
        preimage of $x$ under $g''$. Now set
        $$
            \gamma(x)=\begin{cases}\sup\{f(y):y\in U_x\}&,U_x\neq
            \varnothing\\0&,\text{otherwise}\end{cases}
        $$
        We claim that $\gamma$ is well-defined, that is, the
        supremum in its definition always exists in $X$. For suppose there
        is an $x\in X$ such that the set $\{f(y):y\in U_x\}$ is
        unbounded. Choose a large subset $B\subseteq U_x$
        making the sequence $(f(y):y\in B)$ one-one. Take further
        a large
        $C\subseteq B$ so that all components in the sequence
        $(y:y\in C)$ are either one-one or constant. Set
        $A=\{i\in N: (y_i: i\in C) \text{ is injective}\}$.
        Obviously, $A$ is $f$-wild; therefore it is $g''$-insane.
        Since $g''$ is also monotone, the set $\{g''(y):y\in C\}$
        is large, contradicting the fact that $g''$ is constant
        on $U_x$. Thus, $\gamma$ is well-defined and clearly
        $f\leq h\in\fto{g}$ where $h=\gamma\circ g''$.
    \end{proof}
    \begin{proof}[Proof of Theorem \ref{THM:mainWild}]
        The assertion is an immediate consequence of the preceding
        lemma and the fact that all clones above $\U$ are
        downward closed.
    \end{proof}
    \begin{rem}\label{REM:noConverse}
        Unfortunately, the converse does not hold: If $f,g\in\On$ and $f\in\fto{g}$
        then it need not be true that $f\lw g$. We will see an example at the end of the section.
    \end{rem}
\subsection{$\med_3$ and $T_1$ generate $\pto$}
    We are now ready to prove the explicit description of $\pto$.
    \begin{defn}\label{DEF:Mnk}
        For all $n\geq 1$ and all $1\leq k\leq n$ we define a function
        $$
            m^n_k(x_1,...,x_n)=x_{j_k}\quad ,\text{if} \,\, x_{j_1} \leq ... \leq
            x_{j_n}.
        $$
        For example, $m^n_n$ is the maximum function $\max_n$
        and $m^n_1$ the minimum function $\min_n$ in $n$ variables. Note that
        $\min_n\in Pol(T_1)$ (it is even almost unary) but
        $\max_n\nin\pto$ (and hence $\fto{\max_n}=\OO$).
        If
        $n$ is an odd number then we call $m^n_{\frac{n+1}{2}}$ the
        $n$-th median function and denote this function by $\med_n$.
    \end{defn}
    For fixed odd $n$ it is easily verified (check the wild sets and apply Theorem \ref{COR:FirstptoDescription})
    that $\med_n$ it is the largest of
    the $m^n_k$ which still lies in $Pol(T_1)$: $m^n_k\in
    Pol(T_1)$ iff $k\leq \frac{n+1}{2}$. It is for this reason
    that we are interested in the median functions on our quest
    for a nice generating system of $Pol(T_1)$. As a consequence
    of the following fact from \cite{Pin03} it does not matter which of the median
    functions we consider.
    \begin{fact}
        Let $k,n \geq 3$ be odd natural numbers. Then $\med_k\in \langle \{\med_n\} \rangle$. In
        other words, a clone contains either no median function or
        all median functions.
    \end{fact}

    The following lemma states that within the restrictions of
    functions of $\pto$ (Fact \ref{FAC:goldsternShelah}), we can construct functions of arbitrary
    wildness with the median.
    \begin{lem}\label{LEM:verkettetesSystem}
        Let $n\geq 1$ and let $\A=\{A_1,...,A_k\}\subseteq \P(N)$ be a
        set of subsets of $N$ with the property that $A_i\cap A_j\neq \varnothing$ for all $1\leq i,j\leq k$.
        Then there exists monotone $t_\A\in\clf{\med_3}\un$ such that all
        members of $\A$ are $t_\A$-insane.
    \end{lem}
    \begin{proof}
        We prove this by induction over the size $k$ of $\A$. If
        $\A$ is empty there is nothing to show. If $k=1$, we can
        set $t_\A=\pi^n_i$, where $i$ is an arbitrary element of
        $A_1$. Then $A_1$ is obviously $t_\A$-insane. If $k=2$,
        then define $t_\A=\pi^n_i$, where $i\in A_1\cap A_2$ is
        arbitrary. Clearly, both
        $A_{1}$ and $A_{2}$ are $t_\A$-insane. Finally, assume
        $k\geq 3$. By induction hypothesis, there exist monotone terms
        $t_\B,t_\C,t_\D\in\clf{\med_3}\un$ for the sets $\B=\{A_1,...,A_{k-1}\}$,
        $\C=\{A_1,...,A_{k-2},A_k\}$ and $\C=\{A_{k-1},A_k\}$ such
        that all sets in $\B$ (and $\C,\D$ respectively) are $t_\B$-insane ($t_\C$-insane,
        $t_\D$-insane). Set
        $$
            t_\A=\med_3(t_\B,t_\C,t_\D).
        $$
        Then each $A_i$ is insane for two of the three terms in
        $\med_3$. Thus, if we fix the variables outside $A_{i}$
        to arbitrary values, then at least two of the three subterms
        in
        $\med_3$ are still unbounded and so is $t_\A$ by the
        monotonicity of its subterms.
        Hence, every $A_{i}$ is $t_\A$-insane, $1\leq i\leq k$. Obviously $t_\A$ is
        monotone.
    \end{proof}

    \begin{lem}\label{tf}
        Let $f\in\pton$. Then there exists $t_f\in\clf{\med_3}$
        such that $f\lw t_f$.
    \end{lem}
    \begin{proof}
        Write $\A=\{A_1,...,A_k\}$ for the set of $f$-wild
        subsets of $N$. By Fact \ref{FAC:goldsternShelah}, $A_i\cap A_j\neq \varnothing$ for all
        $1\leq i,j\leq k$. Apply the preceding lemma to $\A$.
    \end{proof}
    \begin{thm}\label{ptoMain}
        $\pto=\fto{\med_3}$.
    \end{thm}
    \begin{proof}
        It is clear that $\pto\supseteq\fto{\med_3}$. On the
        other hand we have just seen that  if $f\in\pto$, then there
        exists $t_f\in\clf{\med_3}$ such that $f\lw t_f$, whence
        $f\in \fto{\med_3}$.
    \end{proof}
    \begin{cor}
        $\pto$ is the $\leq$-downward closure of the clone generated by
        $\med_3$ and the unary functions $\Oo$.
    \end{cor}
    \begin{proof}
        Given $f\in\pto$, by Lemma \ref{tf} there exists $t_f\in\clf{\med_3}$ such that $f\lw t_f$.
        By Lemma \ref{LEM:verkettetesSystem}, $t_f$ is monotone and each $t_f$-wild set is in
        fact even $t_f$-insane. Now one follows the proof of Lemma
        \ref{leq} to obtain $\gamma\in\Oo$ such that $f\leq
        \gamma\circ t_f$.
    \end{proof}
    \begin{cor}
        $\pto=\cl{\{\med_3,p_{\Delta}\}\cup\Oo}$. In particular,
        $\pto$ is finitely generated over the unary functions.
    \end{cor}
    \begin{proof}
        Remember that $\cl{\{p_\Delta\}\cup\Oo}=\cto$ (Fact \ref{FAC:pDelta}) and apply
        Theorem \ref{ptoMain}.
    \end{proof}
        Now we can give the example promised in Remark
        \ref{REM:noConverse}. Set
         $$
            g(x_1,...,x_4)=\med_3(x_1,x_2,x_3)
         $$
         and
         $$
            f(x_1,...,x_4)=\med_5(x_1,x_1,x_2,x_3,x_4).
         $$
         It is obvious that $\fto{g}=\fto{\med_3}=\pto$. Next observe that
         $\fto{f}\subseteq\fto{\med_5}=\pto$ and that
         $f(x_1,x_2,x_3,x_3)=\med_3$ which implies
         $\pto=\fto{\med_3}\subseteq\fto{f}$. Thus,
         $\fto{g}=\fto{f}$.
         Consider on the other hand the 2-element wild sets of the
         two functions: Exactly $\{1,2\},\{1,3\}$ and \{2,3\} are
         $g$-wild, and $\{1,2\},\{1,3\},\{1,4\}$ are the wild sets of two elements for $f$.
         Now the intersection of first group is empty, whereas
         the one of the second group is not; so there is no permutation
         of the set $\{1,2,3,4\}$ which takes the first group to
         the second or the other way. Hence, neither $f\lw g$ nor
         $g\lw f$.
\end{section}

\begin{section}{The interval $[\U,\O]$}
\subsection{A chain in the interval}
    Now we shall show that the open interval $(\cto,\pto)$ is not empty by
    exhibiting a countably infinite descending chain therein with
    intersection $\U$.
    \begin{nota}
        For a natural number $n\geq 2$, we write
        $\M_n=\cl{\{m^n_2\}\cup T_1}$.
    \end{nota}
     Observe that since $m_2^2=\ma_2\nin \pto$, Fact \ref{FAC:poltone} implies that $\M_2=\OO$.
     Moreover, $m_2^3=\med_3$ and hence, $\M_3=\pto$.
    \begin{lem}\label{kLeqn}
        Let $n\geq 2$. Then $\M_n^{(k)}=\U^{(k)}$ for all $1\leq k<n$. That is,
        all functions in $\M_n$ of arity less than $n$ are
        almost unary.
    \end{lem}
    \begin{proof}
        Given $n,k$ we show by induction over terms that if
        $t\in\M_n\uk$, then $t$ is almost unary. To start the
        induction we note that the only $k$-ary functions in the
        generating set of $\M_n$ are almost unary. Now assume
        $t=f(t_1,t_2)$, where $f\in T_1$ and $t_1,t_2\in\M_n\uk$.
        By induction hypothesis, $t_1$ and $t_2$ are almost
        unary and so is $t$ as the almost unary functions are closed under
        composition. Finally, say $t=m_2^n(t_1,...,t_n)$,
        where the $t_i$ are almost unary $k$-ary functions, $1\leq i\leq n$.
        Since $k<n$, there exist $i,j\in N$ with $i\neq j$, $l\in\{1,...,k\}$ and
        $\gamma,\delta\in\Oo$ such that $t_i\leq\gamma(x_l)$ and
        $t_j\leq\delta(x_l)$. Then, $t\leq
        \max(\gamma,\delta)(x_l)$ and so $t$ is almost unary as well.
    \end{proof}
    \begin{cor}\label{MnNinMn1}
        If $n\geq 2$, then $m_2^n\nin\M_{n+1}$. Consequently,
        $\M_n\nsubseteq\M_{n+1}$.
    \end{cor}
    \begin{lem}\label{Mn1InMn}
        If $n\geq 2$, then $m_2^{n+1}\in\M_{n}$. Consequently,
        $\M_{n+1}\subseteq\M_n$.
    \end{lem}
    \begin{proof}
        Set
        $$
            f(x_1,...,x_{n+1})=m_2^n(x_1,...,x_n)\in\M_n.
        $$
        Then every $n$-element subset of $\{1,...,n+1\}$ is
        $f$-wild. Hence, $m_2^{n+1}\lw f$ and so
        $m_2^{n+1}\in\fto{f}\subseteq\M_n$.
    \end{proof}
    \begin{thm}
        The sequence $(\M_n)_{n\geq 2}$ forms a countably infinite descending chain:
        $$
        \OO=\M_2\supsetneqq\M_3=Pol(T_1)\supsetneqq\M_4\supsetneqq
        ...\supsetneqq\M_n\supsetneqq\M_{n+1}\supsetneqq ...
        $$
        Moreover,
        $$
            \bigcap_{n\geq 2} \M_n=\U.
        $$
    \end{thm}
    \begin{proof}
        The first statement follows from Corollary \ref{MnNinMn1}
        and Lemma \ref{Mn1InMn}. The second statement a direct consequence of
        Lemma \ref{kLeqn}.
    \end{proof}

\subsection{Finally, this is the interval}

    We will now prove that there are no more clones in the
    interval $\niceint$ than the ones we already exhibited. We
    first state a technical lemma.

    \begin{lem}
        Let $f\in\On$ be a monotone function such that all
        $f$-wild subsets of $N$ are $f$-insane. Define for $i,j\in N$ with $i\neq j$ functions
        $$
            f^{(i,j)}(x_1,...,x_n)=f(x_1,...,x_{i-1},x_j,x_{i+1},...,x_n)
        $$
        which replace the $i$-th by the $j$-th component and calculate
        $f$. Then the following implications hold for all $f$-wild $A\subseteq N$ and all $i,j\in N$ with $i\neq j$:
        \begin{itemize}
        \item[(i)]{If $i\nin A$, then $A$ is
        $f^{(i,j)}$-insane.}
        \item[(ii)]{If $j\in A$, then $A$ is
        $f^{(i,j)}$-insane.}
        \end{itemize}
    \end{lem}
    \begin{proof}
        We have to show that if we fix the variables outside $A$ to constant values, then $f^{(i,j)}$ is still
        unbounded; because $f$ is monotone, we can assume all values are fixed to
        $0$. Fix a sequence $(\alpha_\xi:\xi\in X)$ of
        elements of $X^n$ such that all components outside $A$ are
        zero for all tuples of the sequence and such that
        $(f(\alpha_\xi):\xi\in X)$ is unbounded. Define a sequence
        of $n$-tuples $(\beta_\xi:\xi\in X)$ by
        $$
            (\beta_\xi)^n_k=\begin{cases}0&,k\nin
            A\\\xi&,\ow\end{cases}
        $$
        For each $\xi\in X$ there exist a $\lambda\in X$ such that
        $\alpha_\xi\leq\beta_\lambda$. Then $f(\alpha_\xi)\leq
        f(\beta_\lambda)$. In either of the cases (i) or
        (ii), $f(\beta_\lambda)\leq f^{(i,j)}(\beta_\lambda)$.
        Thus,
        $(f^{(i,j)}(\beta_\xi):\xi\in X)$ is unbounded.
    \end{proof}

    \begin{lem}
        Let $f\in\On$ not almost unary. Then there exists $n_0\geq
        2$ such that $\fto{f}=\fto{m^{n_0}_2}$.
    \end{lem}
    \begin{proof}
        We shall prove this by induction over the arity $n$ of
        $f$. If $n=1$, there are no not almost unary functions so there is nothing to
        show. Now assume our assertion holds for all $1\leq k<n$.
        We distinguish two cases:

        First, consider $f$ such that
        all $f$-wild subsets of $N$ have size at least $n-1$. Then
        $f\aw m^n_2$ and so $\fto{f}=\fto{m^{n}_2}$.

        Now assume
        there exists an $f$-wild subset of $N$ of size $n-2$, say
        without loss of generality that $\{2,...,n-1\}$ is such a
        set. By Lemma \ref{LEM:MonotoneAndInsane} and Theorem \ref{THM:mainWild} there exists a monotone $\hat{f}$
        with $\fto{f}=\fto{\hat{f}}$ and with the property that all $f$-wild subsets of $N$ are
        $\hat{f}$-insane. Since we could replace $f$ by $\hat{f}$, we assume
        that $f$ is monotone and that all $f$-wild sets are
        $f$-insane.

        Consider the $f^{(i,j)}$ as defined in the preceding lemma. Formally, these functions
        are still $n$-ary, but in fact they depend only on $n-1$ variables. Thus, all
        of the $f^{(i,j)}$ which are not almost unary satisfy the induction
        hypothesis. Set
        $$
            n_0=\min\{k:\exists i,j\in N\,
            \fto{f^{(i,j)}}=\fto{m^{k}_2}\}.
        $$
        The minimum is well-defined: Because $\{2,...,n-1\}$ is $f$-insane,
        $f^{(n,1)}$ is not almost unary so that it generates the same clone
        as some $m^n_2$ modulo $T_1$; thus, the set is
        not empty. Clearly, $m^{n_0}_2\in\fto{f}$. We show that $m_2^{n_0}$ is strong enough to
        generate $f$. Since $\M_n\subseteq\M_{n_0}$ for all $n\geq
        n_0$ we have $f^{(i,j)}\in\fto{m^{n_0}_2}$ for all $i,j\in
        N$ with $i\neq j$. Now define
        $$
            t(x_1,...,x_n)=f^{(n,1)}(x_1,f^{(1,2)},f^{(1,3)},...,f^{(1,n-1)})\in\fto{m^{n_0}_2}.
        $$
        We claim that $f\lw t$. Indeed, let $A\subseteq
        N$ be $f$-wild and whence $f$-insane by our assumption.

        If $1\nin A$, then $A$ is
        $f^{(1,j)}$-insane for all $2\leq j\leq n-1$ by the preceding lemma. So $A$ is
        insane for all components in the definition of $t$ except
        the first one. Hence, because $f$ is monotone, $A$ must be $t$-insane as
        otherwise $f^{(n,1)}$ would be almost unary.

        If $1\in A$, then by the preceding lemma $A$ is still $f^{(1,j)}$-insane whenever $j\in A$.
        Thus, increasing the components with
        index in $A$ increases the first component in $t$ plus all
        subterms $f^{(1,j)}$ with $j\in A$; but by the definition of $f^{(n,1)}$, that is the same
        as increasing the variables $A\cup\{n\}\supseteq
        A$ in $f$. Whence, $A$ is $t$-insane.

        This proves $f\lw t$ and thus $f\in\fto{m^{n_0}_2}$.
    \end{proof}

    So here it is, the interval and the end of our quest.

    \begin{thm}\label{THM:thisIsTheInterval}
        Let $\C\supsetneqq \U$ be a clone. Then there exists
        $n\geq 2$ such that $\C=\M_n$.
    \end{thm}

    \begin{proof}
        Set
        $$
            n_\C=\min\{n\geq 2:\M_n\subseteq\C\}.
        $$
        Since $\C$ contains a function which is not almost unary,
        the preceding lemma implies that
        the set over which we take the minimum is nonempty.
        Obviously, $\M_{n_\C}\subseteq\C$. Now let $f$ be an
        arbitrary function in $\C$ which is not almost unary. Then
        by the preceding lemma, there exists $n_0$ such that
        $\fto{m^{n_0}_2}=\fto{f}$. Clearly, $n_0\geq n_\C$ so that
        $f\in\M_{n_0}\subseteq\M_{n_\C}$.
    \end{proof}

    We state a lemma describing how the $k$-ary
    parts of the $\M_n$ for arbitrary $k$ relate to each other.
    \begin{lem}
         Let $m>n\geq 2$ and $k\geq 2$. If $k\geq n$ (that is, if $\M_n\uk$ is nontrivial), then
         $\M_n\uk\supsetneqq\M_m\uk$.
    \end{lem}
    \begin{proof}
        We know that $\M_n\uk\supseteq\M_m\uk$. To see the
        inequality of the two sets, observe that
        $$
            f(x_1,...,x_k)=m^n_2(x_1,...,x_n)
        $$
        is an element of $\M_n\uk$ but definitely not one of
        $\M_m\uk$.
    \end{proof}

    \begin{cor}
        Let $k\geq 2$. Then
        $$
        \M_2\uk\supsetneqq\M_3\uk\supsetneqq
        ...\supsetneqq\M_k\uk\supsetneqq\M_{k+1}\uk=\U\uk
        $$
        Consequently, there are $k$
        different $k$-ary parts of clones of the interval $[\U,\O]$ for each
        $k$.
    \end{cor}

    In general, if $\C$ is a clone, then
    $$
        Pol(\C^{(1)})\supseteq Pol(\C^{(2)})\supseteq ... \supseteq
        Pol(\C^{(n)})\supseteq ...
    $$
    Moreover,
    $$
        Pol(\C\un)\un=\C\un \quad\text{ and }\quad \bigcap_{n\geq 1} Pol(\C\un)=\C.
    $$
    It is natural to ask whether or not for $\C=\U$ this chain coincides with
    the chain we discovered.
    \begin{thm}
        Let $n\geq 1$. Then $\M_{n+1}=Pol(\U\un)$.
    \end{thm}
    \begin{proof}
        Clearly, $\M_2=Pol(\U^{(1)})=\O$, so assume $n\geq 2$.
        Consider $m^{n+1}_2$ and let $f_1,...,f_{n+1}$ be
        functions in $\U\un$. Then two of the $f_j$ are bounded by
        unary functions of the same variable. Thus
        $m^{n+1}_2(f_1,...,f_{n+1})$ is bounded by a unary function of
        this variable. This shows $m^{n+1}_2\in Pol(\U\un)$ and
        hence $\M_{n+1}\subseteq Pol(\U\un)$. Now consider $m^n_2$ and observe that
        $m^n_2\nin\U\un=Pol(\U\un)\un$; this proves $\M_{n}\nsubseteq
        Pol(\U\un)$. Whence, $\M_{n+1}=Pol(\U\un)$.
    \end{proof}
\subsection{The $m^n_k$ in the chain}
    As an example, we will show where the clones generated by the $m^n_k$ (as in Definition \ref{DEF:Mnk}) and $T_1$ can be found in
    the chain.
    \begin{nota}
        For $1\leq k\leq n$ we set $\M_n^k=\fto{m_k^n}$.
    \end{nota}
    Note that if $k=1$, then $\M_n^k=\U$, and if $k>\frac{n+1}{2}$,
    then
    $\M_n^k=\OO$. Observe also that $\M_n=\M^2_n$ for all $n\geq 2$.
    \begin{nota}
        For a positive rational number $q$ we write
        $$
          \lfloor q \rfloor = \max\{n\in\mathbb{N}:\,n\leq q\}
        $$
        and
        $$
            \lceil q \rceil = \min\{n\in\mathbb{N}:\,q\leq n\}.
        $$
        The remainder of the division $\frac{n}{k}$ we denote by the
        symbol $R(\frac{n}{k})$.
    \end{nota}
    \begin{lem}
        Let $2\leq k\leq \frac{n+1}{2}$ and let $t\in\M_n^k$ not almost unary. Then all $t$-wild
        subsets of $N_t$ have size at least
        $\frac{n}{k-1}-1$.
    \end{lem}
    \begin{proof}
        Our proof will be by induction over terms. If $t=m_k^n$,
        then all $t$-wild subsets of $N$ have at least $n-k+1$
        elements in accordance with our assertion. For the
        induction step, assume $t=f(t_1,t_2)$, where $f\in T_1$,
        say $f(x_1,x_2)\leq\gamma(x_1)$ for some $\gamma\in\Oo$.
        Then $t$ inherits the asserted property from $t_1$.
        Finally we consider the case where $t=m_k^n(t_1,...,t_n)$.
        Suppose towards contradiction there exists $A\subseteq
        N_t$ $t$-wild with $|A|<\frac{n}{k-1}-1$.
        There have to be at least $n-k+1$ terms $t_j$ for which
        $A$ is $t_j$-wild so that $A$ can be $t$-wild. By
        induction hypothesis, these $n-k+1$ terms are almost
        unary and bounded by a unary function of a variable with index in $A$.
        From the bound on the size of $A$ we conclude that there
        must be an index in $A$ so that at least
        $$
            \lceil\frac{n-k+1}{|A|}\rceil>\lceil\frac{n-k+1}{\frac{n}{k-1}-1}\rceil=k-1
        $$
        of the terms $t_j$ are bounded by an unary function of the
        same variable. But if $k$ of the $t_j$ have the same one-element strong set,
        then $t$ is bounded by a unary function of this variable as
        well in contradiction to the assumption that $t$ is not
        almost unary.
    \end{proof}
    \begin{cor}\label{COR:Mnknsubseteq}
        Let $2\leq k\leq \frac{n+1}{2}$. Then $\M_{\lceil\frac{n}{k-1}\rceil-1}\nsubseteq\M_n^k$.
    \end{cor}
    \begin{proof}
        With the preceding lemma it is enough to observe that
        $m_2^{\lceil\frac{n}{k-1}\rceil-1}\in\M_{\lceil\frac{n}{k-1}\rceil-1}$
        has a wild set of size $\lceil\frac{n}{k-1}\rceil-2$.
    \end{proof}
    So we identify now the $\M_j$ which $\M_n^k$ is
    equal to.
    \begin{lem}\label{LEM:M_quotientFromMnk}
        Let $2\leq k\leq n$. Then $\M_{\lceil\frac{n}{k-1}\rceil}\subseteq\M_n^k$.
    \end{lem}
    \begin{proof}
        It suffices to show that $m_k^n$ generates
        $m_2^{\lceil\frac{n}{k-1}\rceil}$. But this is easy:
        $$
            m_2^{\lceil\frac{n}{k-1}\rceil}=m_k^n(x_1,...,x_1,x_2,...,x_2,...,x_{\lceil\frac{n}{k-1}\rceil},...,x_{\lceil\frac{n}{k-1}\rceil}),
        $$
        where $x_j$ occurs $k-1$ times if $1\leq j\leq
        \lfloor\frac{n}{k-1}\rfloor$ and $R(\frac{n}{k-1})<k-1$ times
        if $j=\lfloor\frac{n}{k-1}\rfloor+1$. For if we evaluate the function for a $\lceil\frac{n}{k-1}\rceil$-tuple with
        $x_{j_1}\leq
        ...\leq x_{j_{\lceil\frac{n}{k-1}\rceil}}$, then $x_{j_1}$
        occurs at most $k-1$ times in the tuple, but $x_{j_1}$
        together with $x_{j_2}$ occur more than $k$ times; thus,
        the $k$-th smallest element in the tuple is $x_{j_2}$ and
        $m_k^n$ returns $x_{j_2}$.
    \end{proof}
    \begin{thm}
        $\M_n^k=\M_{\lceil\frac{n}{k-1}\rceil}$ for all $2\leq k\leq
        n$.
    \end{thm}
    \begin{proof}
         By Theorem \ref{THM:thisIsTheInterval}, $\M_n^k$ has to be somewhere in the chain $(\M_n)_{n\geq 2}$.
         Because of Corollary \ref{COR:Mnknsubseteq} and Lemma
         \ref{LEM:M_quotientFromMnk} the assertion follows.
    \end{proof}
\subsection{Further on the chain}
    We conclude by giving one simple guideline for where to search
    the clone $\fto{f}$ in the chain for arbitrary $f\in \O$.

    \begin{lem}\label{kElWildSet}
        Let $2\leq k\leq n$ and let $f\in\On$ be a not almost unary function
        which has a $k$-element $f$-wild subset of $N$. Then
        $\M_{k+1}\subseteq\fto{f}$.
    \end{lem}
    \begin{proof}
        We can assume that $\{1,...,k\}$ and all $A\subseteq N$ with $|A|=n-1$ are $f$-insane
        and that $f$ is monotone. Define
        $$
            g(x_1,...,x_{k+1})=f(x_1,...,x_k,x_{k+1},...,x_{k+1})\in\fto{f}.
        $$
        Let $A\subseteq\{1,...,k+1\}$ with $|A|=k$ be given. If
        $A=\{1,...,k\}$ then $A$ is $f$-wild and so it is
        $g$-wild. Otherwise $A$ contains $k+1$ and so it affects $n-1$ components in the
        definition of $g$. Therefore $A$ is
        $g$-wild by Lemma \ref{LEM:almostUnaryCriterion}. Hence, $m_2^{k+1}\lw g$
        and so $\M_{k+1}\subseteq\fto{g}\subseteq\fto{f}$.
    \end{proof}
    \begin{rem}
        Certainly it is not true that if the smallest wild set of
        a function $f\in\O$ has $k$ elements, then
        $\M_{k+1}=\fto{f}$. The $m^n_k$ are an example.
    \end{rem}
    \begin{cor}\label{COR:2elwildset}
        Let $f\in\pto$ not almost unary and such that there exists a
        2-element $f$-wild subset of $N$. Then $\fto{f}=\pto$.
    \end{cor}
\subsection{A nice picture}
    If $X$ is countably infinite or weakly compact, we can draw the situation we ran into like this.
\newpage
\vspace*{4mm}

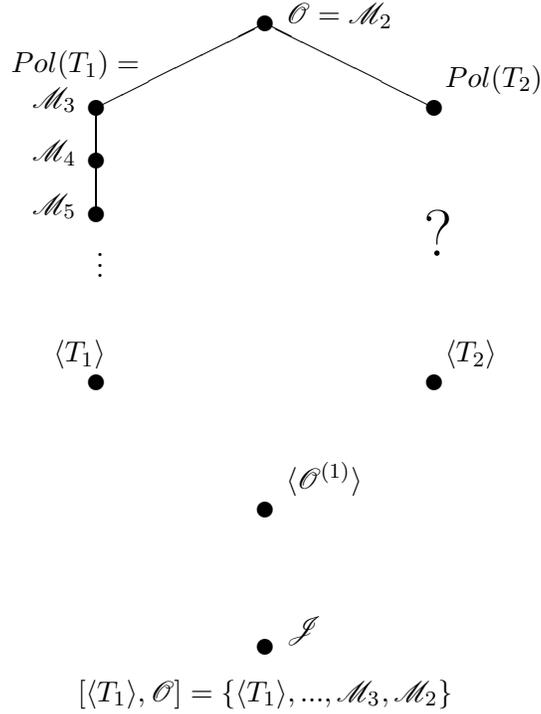
\setlength{\unitlength}{0.8\unitlength}
    \begin{center}
    \begin{picture}(200,300)
    \put(100,5){\circle*{8}}
    \put(110,10){$\J$}

    \put(100,70){\circle*{8}}
    \put(110,80){$\cl{\Oo}$}

    \put(100,300){\circle*{8}}
    \put(110,300){$\OO=\M_2$}

    \put(180,260){\circle*{8}}
    \put(185,270){$Pol(T_2)$}

    \put(20,130){\circle*{8}}
    \put(0,140){$\cto$}

    \put(180,130){\circle*{8}}
    \put(185,140){$\langle T_2\rangle$}

    \put(20,260){\line(2,1){80}}
    \put(180,260){\line(-2,1){80}}

    \put(175,190){\Huge{?}}

    \put(20,260){\circle*{8}}
    \put(-20,277){$\pto=$}
    \put(-10,260){$\M_3$}

    \put(20,210){\line(0,1){50}}

    \put(20,235){\circle*{8}}
    \put(-10,235){$\M_4$}

    \put(20,210){\circle*{8}}
    \put(-10,210){$\M_5$}

    \put(20,180){$\vdots$}

    \end{picture}
    \end{center}
    $$[\cto,\O]=\{\cto,...,\M_3,\M_2\}$$

\end{section}
\setlength{\unitlength}{1.25\unitlength}


\begin{thebibliography}{10}

\bibitem{DR85}
R.~O. Davies and I.~G. Rosenberg, \emph{Precomplete classes of
operations on an
  uncountable set}, Colloq. Math. \textbf{50} (1985), 1--12.

\bibitem{Gav65}
G.~P. Gavrilov, \emph{On functional completeness in
countable-valued logic
  ({Russian})}, Problemy Kibernetiki \textbf{15} (1965), 5--64.

\bibitem{Gol02}
M.~Goldstern, \emph{Analytic clones}, preprint.

\bibitem{GS03}
M.~Goldstern and S.~Shelah, \emph{Clones from creatures},
preprint.

\bibitem{GS022}
\bysame, \emph{Clones on regular cardinals}, Fundamenta
Mathematicae
  \textbf{173} (2002).

\bibitem{Kec95}
A.~Kechris, \emph{Classical descriptive set theory}, Springer,
1995.

\bibitem{Pin03}
M.~Pinsker, \emph{The clone generated by median functions},
preprint.

\bibitem{Pin02}
\bysame, \emph{Rosenberg's characterization of maximal clones},
Master's
  thesis, Vienna University of Technology, 2002.

\bibitem{Ros70}
I.~G. Rosenberg, \emph{\"{U}ber die funktionale
{Vollst\"{a}ndigkeit} in den
  mehrwertigen {Logiken}}, Rozpravy \v{C}eskoslovensk\'{e} Akad. v\v{e}d, Ser.
  Math. Nat. Sci. \textbf{80} (1970), 3--93.

\bibitem{Ros76}
\bysame, \emph{The set of maximal closed classes of operations on
an infinite
  set ${A}$ has cardinality $2^{2^{|A|}}$}, Arch. Math. (Basel) \textbf{27}
  (1976), 561--568.

\end{thebibliography}
\end{document}